\documentclass[12pt]{article}
\usepackage{amsxtra,amssymb,amsthm,amsmath,latexsym}

\theoremstyle{plain}
\newtheorem{theorem}{Theorem}

\newtheorem{lemma}{Lemma}
\newtheorem{remark}{Remark}

\newcommand{\refT}[1]{Theorem~\ref{T:#1}}
\newcommand{\refS}[1]{Section~\ref{S:#1}}

\newcommand{\refL}[1]{Lemma~\ref{L:#1}}

\def\d{\delta}
\def\dota{\dot{a}}
\def\dotmu{\dot{\mu}}

\def\dotg{\dot{g}}
\def\dotu{\dot{u}}
\def\dotw{\dot{w}}
\def\tildeA{\widetilde A}

\def\tildeT{\widetilde T}

\def\oH{{\overset{\circ}{H}}}
\def\oH1{{\overset{\circ}{H}\kern-.02in{}^1}}

\def\loc{{\hbox{\,loc\,}}}
\def\bee{\begin{equation*}}
\def\eee{\end{equation*}}
\def\be{\begin{equation}}
\def\ee{\end{equation}}

\begin{document}
\title{ Dynamical systems method (DSM) for general nonlinear 
equations.}

\author{A.G. Ramm\\
 Mathematics Department, Kansas State University, \\
 Manhattan, KS 66506-2602, USA\\
ramm@math.ksu.edu}

\date{}
\maketitle\thispagestyle{empty}

\begin{abstract}
\footnote{key words: nonlinear operator equations, iterative 
methods, DSM-dynamical systems method  }

If $F:H\to H$ is a map in a Hilbert space $H$, $F\in C^2_{loc}$,
and there exists $y$, such that $F(y)=0$, $F'(y)\not= 0$,
then equation $F(u)=0$ can be solved by a DSM (dynamical systems method).
This method yields also a convergent iterative method for finding $y$,
and this method converges at the rate of a geometric series. It is not 
assumed
that $y$ is the only solution to $F(u)=0$.
Stable approximation to a solution of the equation $F(u)=f$ is
constructed by a DSM when $f$ is unknown but $f_\d$ is known, where
$||f_\d-f||\leq \d$. 
\end{abstract}

\section{Introduction}\label{S:1}
In this paper a method for solving fairly general class of nonlinear 
operator equations $F(u)=0$ in a Hilbert space is proposed, its 
convergence is proved,
and an iterative method for solving the above equation is constructed.
Convergence of the iterative method is also proved. These results are 
based
on the following assumptions: a) the above equation has a solution $y$, 
possibly non-unique,
b) $F\in C^2_{\loc}$, 
and
c) $F'(y)\not= 0$. No restrictions on the rate of growth of nonlinearity 
are made. The literature on methods for solving nonlinear equations is
large (see, e.g., \cite {OR} and references therein). Most of the results
obtained so far are based on Newton-type methods and their modifications.
There is also a well-developed theory for equations with monotone 
operators (\cite{Z}).
The method used in this paper is a version of the dynamical systems method
which is studied in \cite{R1}. The idea of this method is described 
briefly below.

Let $F:H\to H$ be a map in a Hilbert space. Assume that equation
\be\label{e1} F(u)=0 \ee
has a solution $y$, possibly non-unique, and
\be\label{e2} F'(y)\not= 0, \ee
where $F'$ is the Fr\`echet derivative of $F$. This ssumption means that 
$F'(y)$ is not equal to zero identically on $H$. Thus, it is a weak 
assumption.
Assume that $F\in C^2_{loc}$, i.e.,
\be\label{e3}
\sup_{u\in B(u_0,R)} \|F^{(j)}(u)\| \leq M_j (R) \qquad 0\leq j\leq 2, 
\ee
where $u_0\in H$ is a given element, $R>0$, and no restrictions on the
growth of $M_j(R)$ as $R$ grows are made.
This means that the nonlinearity $F$ can grow arbitrarily fast as
$\|u-u_0\|$ grows. It is known that under these assumptions
equation \eqref{e1} may have no solutions. Thus, we have assumed
that a solution $y$ to \eqref{e1} exists.

We do not assume that $F'(u)$ has a bounded inverse operator.
Therefore the standard Newton-type methods are not applicable. Dynamical
system method (DSM) consists of finding an operator $\Phi$ such that
the problem
\be\label{e4} \dotu=\Phi(t,u), \qquad u(0)=u_0 \ee
has a unique global solution $u(t)$, there exists $u(\infty)$,
and $F(u(\infty))=0$. To ensure the unique local solvability of \eqref{e4}
we assume that 
$$\|\Phi(t,u)-\Phi(t,v)\|\leq L(R) \|u-v\|\quad
\forall u,v\in B(u_0,R).$$
 Then the global existence of the unique local
solution holds if $\sup_t \|u(t)\|<\infty$.

The results of this paper can be summarized in two theorems. Let us denote
\be\label{e5}
A:=F'(u(t)),\,\, T:=A^\ast A, \quad T_a:=T+aI; \quad \tildeA:=F'(y),
\quad \tildeT=\tildeA^\ast \tildeA. \ee
Assume that $a(t)$ is a positive monotonically decaying function,
\be\label{e6}
a(t)>0, \quad \lim_{t\to\infty} a(t)=0,
\quad \frac{|\dota|}{a}\leq\frac{1}{2}. \ee

\begin{theorem}\label{T:1}
If a solution $y$ to equation \eqref{e1} exists, possibly is non-unique,
and assumptions \eqref{e2} and \eqref{e3} hold, 
then $y=u(\infty)$, where 
$u(t)$ solves the following DSM problem:
\be\label{e7}
\dotu =- T^{-1}_{a(t)} [A^\ast F(u(t))+a(t) (u(t)-z)],
\quad u(0)=u_0, \ee
and where $z$ and $u_0$ are suitably chosen.
\end{theorem}

\begin{theorem}\label{T:2}
Under the assumptions of {\rm \refT{1}}, the iterative process
\be\label{e8}
u_{n+1} =u_n-h_n T^{-1}_{a_n} [A^\ast(u_n) F(u_n)+a_n(u_n-z)],
\quad u_0=u_0, \ee
where $h_n>0$ and $a_n>0$ are suitably chosen, generates the
sequence $u_n$ converging to $y$.
\end{theorem}

\begin{remark}\label{R:1}
The suitable choices of $a_n$ and $h_n$ are discussed in the proof of 
{\rm \refT{2}}.
\end{remark}

\begin{remark}\label{R:2}
Essentially, {\rm \refT{1}} says that {\rm any} solvable operator equation
with $C^2_{loc}$ operator, satisfying only a weak assumption
\eqref{e2}, can be solved by a DSM. Condition \eqref{e2} means
that the range of the linear operator $F'(y)$ contains at least one
non-zero element. It allows $F'(y)$ to have infinite-dimensional
null-space.
\end{remark}

In \refS{2} we prove \refT{1} and \refT{2}. In their proofs we use
the following lemmas.

\begin{lemma}\label{L:1}
Assume that $g(t)\geq 0$ is a $C^1([0,\infty)]$ function satisfying
the inequality
\be\label{e9}
\dotg \leq -\gamma(t) g+\alpha(t)g^2+\beta(t),
\quad t\geq 0, \quad \dotg:=\frac{dg}{dt}, \ee
where $\gamma$, $\alpha$ and $\beta$ are nonnegative continuous functions.
Assume that there exists $\mu\in C^1([0,\infty))$, $\mu>0$,
$\lim_{t\to\infty}\mu(t)=\infty$, such that
\be\label{e10}
\hbox{i)\ } \alpha(t)\leq \frac{\mu(t)}{2}
        \left( \gamma(t)-\frac{\dotmu(t)}{\mu(t)}\right),
        \quad
\hbox{ii)\ } \beta\leq\frac{1}{2\mu(t)}
        \left( \gamma(t)-\frac{\dotmu(t)}{\mu(t)}\right),
        \quad
\hbox{iii)\ } g(0)\mu(0)<1. \ee
Then any solution to \eqref{e9} exists on $[0,\infty)$ and
\be\label{e11}
0\leq g(t)<\frac{1}{\mu(t)}, \qquad t\in[0,\infty). \ee
\end{lemma}

\refL{1} is proved in \cite[pp.66-70]{R1}.

\begin{lemma}\label{L:2}
Let $g_{n+1}\leq \gamma g_n+pg^2_n$, $g_0:=m>0$, $0<\gamma<1$, $p>0$.
If $m<\frac{q-\gamma}{p}$, where $\gamma<q<1$, then
$\lim_{n\to\infty} g_n=0$, and $g_n\leq g_0 q^n$.
\end{lemma}

{\it Proof of \refL{2}}. Estimate $g_1\leq \gamma m+pm^2\leq 
qm$
holds if $m\leq \frac{q-\gamma} p$, $\gamma< q <1$. Assume that $g_n\leq 
g_0 q^n$. Then 
$$g_{n+1}\leq \gamma g_0q^n+ 
p(g_0q^n)^2=g_0q^n(\gamma+pg_0q^n)<g_0q^{n+1},
$$
because $\gamma+pg_0q^n<\gamma+pg_0q\leq q$. Lemma 2 is proved. \hfill 
$\Box$

\section{Proofs}\label{S:2}

\begin{proof}[Proof of \refT{1}]
If $F'(y):=A$ is linear and $A\not= 0$, then there exists $v_1\not= 0$,
$v_1=\tildeT v$. By the linearity of $\tildeT,$ every element 
$cv_1$ belongs to the range of $\tildeT$ for any constant 
$c$, because $cv_1=\tildeT(cv)$.
Therefore there exists a $z$ such that $y-z=\tildeT v$,
where $\|v\|>0$ can be chosen arbitrarily small. How small $\|v\|$
should be chosen will become clear later. Let $u(t)-y:=w(t)$, 
$\|w(t)\|:=g(t)$.
Write equation \eqref{e7} as
\be\label{e12}
\dotw=-T^{-1}_{a(t)}\  [A^\ast(F(u)-F(y)) +a(t) w+a(t)(y-z)], \ee
and use the formula $F(u)-F(y)=Aw+K$, where $\|K\|\leq\frac{M_2g^2}{2}$.
Then
\be\label{e13}
\dotw=-w-T^{-1}_{a(t)}\ A^\ast K-a(t) T^{-1}_{a(t)} \tildeT v. \ee
Multiply this equation by $w$ in $H$ and use the estimate
$\|T^{-1}_a A^\ast\|\leq\frac{1}{2\sqrt{a}}$, $a>0$, to get
\bee
g\dotg\leq -g^2+\frac{1}{2\sqrt{a(t)}} \frac{M_2g^2}{2}
  +a(t) \ \|\left( T^{-1}_{a(t)} - \tildeT^{-1}_{a(t)}
  + \tildeT^{-1}_{a(t)}\right)\tildeT\| \|v\|.
\eee
If $a>0$ then  $\|\tildeT^{-1}_{a} \tildeT\|\leq 1$,
\bee
a \|T^{-1}_{a}\| \leq 1, \quad
a\|(T^{-1}_a-\tildeT^{-1}_a)\tildeT\|
   = a\|T^{-1}_a (A^\ast A-\tildeA^\ast \tildeA) \tildeT^{-1}_a\tildeT\|
\leq 2M_1M_2g. \eee
Collecting the above estimates and choosing
$\|v\|$ so that $2M_1M_2\|v\|\leq\frac{1}{2}$,
we obtain
\be\label{e14}
\dotg\leq -\frac{g}{2} + \frac{c_0 g^2}{\sqrt{a(t)}} + a(t) \|v\|,
\qquad c_0:=\frac{M_2}{4}. \ee
Apply \refL{1} to \eqref{e14}. Here $\gamma=\frac{1}{2}$,
$\alpha=\frac{c_0}{\sqrt{a(t)}}$, $\beta=a(t)\|v\|$.
Let $\mu(t)=\frac{\lambda}{\sqrt{a(t)}}$, $\lambda=const>0$.
%
Condition $i)$ of \refL{1} holds if
$\frac{c_0}{\sqrt{a(t)}} \leq\frac{\lambda}{2\sqrt{a(t)}}
  \left(\frac{1}{2}-\frac{1}{2}\frac{|\dota|}{a}\right)$.
This inequality holds if $\lambda\geq 8c_0$,
see the last assumption \eqref{e6}.
Condition $iii)$ holds if $g(0)\frac{\lambda}{\sqrt{a(0)}}<1$.
This inequality holds (for any initial value $g(0)=\|u_0-y\|$)
if $a(0)$ is sufficiently large.
Condition $ii)$ holds if $\sqrt{a(t)}\|v\|\leq\frac{1}{8\lambda}$,
where we have used the last assumption \eqref{e6} again.
This inequality holds if 
$(\ast)\quad 8\lambda\sqrt{a(0)}\|v\|\leq 1$. (Recall that $a(0)\geq a(t)$
due to monotonicity of $a(t)$.)
Inequality $(\ast)$ holds if $\|v\|$ is sufficiently small.
Thus, if $\|v\|$ is sufficiently small, then, by \refL{1}, we get
$g(t)<\frac{\sqrt{a(t)}}{\lambda}$,
so $\|u(t)-y\|\leq \frac{\sqrt{a(t)}}{\lambda}\to 0$ as $t\to\infty$.

\refT{1} is proved.
\end{proof}

\begin{proof}[Proof of \refT{2}]

Let $w_n:=u_n-y$, $g_n:=\|w_n\|$. As in the proof of Theorem 1, we assume
$2M_1M_2\|v\|\leq\frac{1}{2}$ and rewrite (8) as
\bee
w_{n+1}=w_n-h_n\ T^{-1}_{a_n}\ [A^\ast(u_n)(F(u_n)-F(y))+
a_nw_n + a_n(y-z)], \quad w_0=\|u_0-y\|. \eee
Using the Taylor formula
$$F(u_n)-F(y)=A(u_n) w_n+K(w_n), \quad
\|K\|\leq\frac{M_2g_n^2}{2},$$
 the estimate
$\|T^{-1}_{a_n} A^\ast(u_n)\| \leq\frac{1}{2\sqrt{a_n}}$, and the formula
$y-z=\tildeT v$, we get
\be\label{e15}
w_{n+1}=(1-h_n) w_n-h_n T^{-1}_{a_n}
A^\ast(u_n)K(w_n) -h_n a_n T^{-1}_{a_n} \tildeT v. \ee
Taking into account that $\|\tildeT^{-1}_a\tildeT\|\leq 1$, and
$a\|T^{-1}_a\|\leq 1$ if $a>0$,
we obtain
\bee
\|T^{-1}_{a_n} \tildeT v\|
\leq \|(T^{-1}_{a_n}-\tildeT^{-1}_{a_n}) \tildeT \| \|v\|+\|v\|,\eee
and
\bee
\|(T^{-1}_{a_n}-\tildeT^{-1}_{a_n}) \tildeT\|
=\|T^{-1}_{a_n}(\tildeT_{a_n}-T_{a_n}) \tildeT^{-1}_{a_n} \tildeT\|
\leq \frac{2M_1M_2g_n}{a_n}:=\frac{c_1g_n}{a_n}. \eee
Let $c_0:=\frac{M_2}{4}$.
Then we obtain from \eqref{e15} the following inequality:
\bee
   g_{n+1}\leq (1-h_n) g_n+\frac{c_0h_n g^2_n}{\sqrt{a_n}}
   + c_1 h_n \|v\| g_n + h_na_n \|v\|. \eee
We have assumed in the proof of Theorem 1 that 
$c_1\|v\|\leq\frac{1}{2}$. Thus
$$g_{n+1} \leq(1-\frac{h_n}{2}) g_n
+\frac{c_0h_n}{\sqrt{a_n}} g^2_n +h_n a_n \|v\|.$$
Choose $a_n=16 c^2_0 g^2_n$.
Then  $\frac{c_0 g_n}{\sqrt{a_n}}=\frac 1 4$, and
\be\label{e16}
g_{n+1}\leq (1-\frac{h_n}{4}) g_n+16 c_0h_n \|v\| g^2_n,
\qquad g_0=\|u_0-y\|\leq R, \ee
where $R>0$ is defined in \eqref{e3}.
Take $h_n=h\in(0,1)$ and choose $g_0:=m$ such that
$m<\frac{q+h-1}{16c_0 h\|v\|}$, where $q\in(0,1)$ and $q+h>1$.
Then \refL{2} implies 
$$\|u_n-y\|\leq g_0 q^n\to 0 \hbox { as }n\to\infty.$$

\refT{2} is proved.

\end{proof}

\section{Stability of the solution}\label{S:3}
Assume that $F(y)=f$, where the exact data $f$ are not known but
the noisy data $f_\d$ are given, $||f_\d-f||\leq \d$. Then the DSM
yields a stable approximation of the solution $y$ if the stopping time 
$t_\d$ is properly chosen. The DSM is similar to (7):
\be\label{e17}
\dotu_\d =- T^{-1}_{a(t)} [A^\ast (F(u_\d(t))-f_\d)+a(t) (u_\d(t)-z)],
\quad u_\d(0)=u_0, \ee
Let $w_\d:=u_\d(t)-y$, $g_\d(t):=||w_\d||$. As in the proof of Theorem 1
we derive the inequality similar to (14):
\be\label{e18}
\dotg_\d\leq -\frac{g_\d}{2} + \frac{c_0 g_\d^2}{\sqrt{a(t)}} + a(t) 
\|v\|+\frac {\d}{2\sqrt{a(t)}},
\qquad c_0:=\frac{M_2}{4}, \ee
and apply Lemma 1. The only 
difference is in checking condition $ii)$ of Lemma 1. This condition now 
takes the form: 
$$a(t)\|v\|+\frac {\d}{2\sqrt{a(t)}}\leq \frac 
{\sqrt{a(t)}}{8\lambda}.$$
This condition can only be satisfied for $t\in [0,t_\d]$, where 
$t_\d<\infty$. The stopping time $t_\d$
can be determined, for example, from the equation $4\lambda \frac 
{\d}{a(t)}=\frac 1 2$, provided that $v$ is chosen sufficiently small,
so that $8\lambda \sqrt{a(0)}||v||\leq \frac 1 2$, and 
$\lambda \geq 8c_0$ as in the proof of Theorem 1. Then, by Lemma 1, we 
have
$g_\d(t_\d)<\frac  {\sqrt{a(t)}}{\lambda}\to 0$ as $\d\to 0$.
Let us formulate the result.
\begin{theorem}\label{T:3}
Let $u_\d:=u_\d(t_\d)$, where $u_\d(t)$ solves problem (17) and $t_\d$ is 
chosen as above. Then $\lim_{\d\to 0}||u_\d-y||=0$.
\end{theorem}

\end{document}